\documentclass[graybox]{svmult}

\usepackage{type1cm}        
\usepackage{makeidx}         
\usepackage{graphicx}        
\usepackage{multicol}        
\usepackage[bottom]{footmisc}

\usepackage{newtxtext}       %
\usepackage{newtxmath}       

\newcommand{\T}{\mathbb T}
\newcommand{\R}{\mathbb R}
\newcommand{\Z}{\mathbb Z}
\newcommand{\N}{\mathbb N}
\newcommand{\p}{\partial}

\makeindex             

\begin{document}

\title*{One-dimensional turbulence with Burgers}
\author{Roberta Bianchini and Anne-Laure Dalibard}
\institute{Roberta Bianchini \at  Sorbonne-Universit\'e, CNRS, Universit\'e de Paris, Laboratoire Jacques-Louis Lions (LJLL), F-75005 Paris, France \& Consiglio Nazionale delle Ricerche, IAC, via dei Taurini 19, I-00185 Rome (Italy)  \email{r.bianchini@iac.cnr.it}
\and Anne-Laure Dalibard \at    Sorbonne-Universit\'e, CNRS, Universit\'e de Paris, Laboratoire Jacques-Louis Lions (LJLL), F-75005 Paris, France \email{dalibard@ljll.math.upmc.fr}}
%
%
\maketitle

\abstract{Gathering together some existing results, we show that the solutions to the one-dimensional Burgers equation converge for long times towards the stationary solutions to the steady Burgers equation, whose Fourier spectrum is not integrable. This is one of the main features of wave turbulence.}

\section{Introduction}

In this paper, we are interested in the long time behavior of the forced inviscid Burgers  equation. 
Indeed, as advertised by Menzaque, Rosales, Tabak and Turner, this equation is the simplest good example simulating the mechanism of weak wave turbulence \cite{MRTT}: ``\emph{the nonlinear term in \eqref{Burgers} has two combined functions: to transfer energy among the various Fourier components of $u$ and to dissipate energy at shocks. Thus the inertial cascade and the dissipation are modeled by the same term}.''

In \cite{Laure}, Yves Colin de Verdi\`ere and Laure Saint-Raymond are able to show that the \emph{linear} and inviscid internal waves model in a trapezoidal domain with monochromatic forcing has continuous spectrum. As a consequence, the behavior of the system for long times is different from the integrable case, because of the quasi-resonance mechanism. In particular, the solutions converge to the generalized eigenfunctions of the linear operator, which are rough functions contained in some negative order Sobolev space. Passing to the Fourier representation, their spectrum is therefore not integrable, like the case of the so-called \emph{Kolmogorov solutions} (due to Zakharov) to the weak wave turbulence model, see for instance \cite{Majda, Nazarenko}.\\
In this note, we will rather discuss a one-dimensional \emph{nonlinear} equation in a regular domain (a finite interval of the real line with periodic boundary conditions), whose solutions converge for long times to some steady solutions which are not integrable from the Fourier side.\\\\

\section{Long-time behavior of the forced Burgers equation with steady forcing}

We are interested in the long time behavior of  the solutions of the forced Burgers equation
\begin{equation}\label{Burgers}
\partial_t u + u\partial_x u = f(x),\quad x\in \T, \ t>0,
\end{equation}
where $\T:=\R/(2\pi \Z)$,
for some periodic potential $f$ with zero average. 
In this note, to illustrate the analysis, we will mostly work with $f(x)=-\kappa_0 \cos(\kappa_0 x/2) \sin(\kappa_0 x/2) $, for some integer $\kappa_0\in \mathbb N$, but we will also make comments in the case of an arbitrary potential $f$.

We endow \eqref{Burgers} with an initial data $u(t=0)=u_0\in L^\infty(\T)$.
It is classical that the average of the solution $u$ of \eqref{Burgers} is preserved by the evolution.
We set $p:=\langle u_0\rangle=\langle u(t)\rangle$, where $\langle g\rangle:= \frac{1}{2\pi}\int_0^{2\pi} g$ for $g\in L^1(\T)$.
Setting $v=u-p$, we can re-write \eqref{Burgers} either as
\[
\p_t v + \frac{1}{2}\p_x (p+ v)^2=f(x),
\]
or, following \cite{MRTT}, as
\begin{equation}
   \p_t v + p\p_x v + v\p_xv = f(x), 
\end{equation}\label{Burgers-bis}
and $v$ now has zero average. The first formulation will be useful when we use the equivalence with Hamilton-Jacobi equation, while the second one will be relevant when we investigate resonance mechanisms.


We give in this note a proof of the following result:
\begin{theorem}
Let $u_0\in L^\infty(\T)$, and let $p=\langle u_0\rangle$. Let $u\in L^\infty_\text{loc}(\R_+, L^\infty(\T))$ be the unique entropy solution of \eqref{Burgers} such that $u_{|t=0}=u_0$.

Then there exists $\bar u\in L^\infty(\T)$ a stationary entropy solution of \eqref{Burgers} such that
\[
u(t)\to \bar u\quad \text{as }t\to \infty\text{ in } L^q(\T),\ 1\leq q <+\infty.
\]
Furthermore, there exists $p_{cr}\geq 0$ ($p_{cr}=2\sqrt{2}/\pi$ in the case when $f(x)=-\kappa_0 \cos(\kappa_0 x/2) \sin(\kappa_0 x/2) $) such that $\bar u \in \mathcal C(\T)$ if $|p|\geq p_{cr}$, and $\bar u$ is discontinuous if $|p|<p_{cr}$.
\label{thm-Burgers}
\end{theorem}

We do not claim that the above result is new, and in fact, we will rely on previous results (especially concerning the long-time behavior or the homogenization of Hamilton-Jacobi equations) to prove Theorem \ref{thm-Burgers}. Following \cite{MRTT}, we will also study the case when the forcing $f$ is of the form $f(x-\omega t).$ 

In order to characterize the long-time behavior of $u$  in $L^1(\T)$, we develop here an argument which glues together several  existing  results from the literature:
\begin{itemize}
\item First, because of the non-degeneracy of the flux (and in fact, in this case, of its convexity), any sequence $u(t_n+\cdot)$ is compact in $\mathcal C([0,T],L^1(\T))$  for all $T>0$. 
This allows us to consider the $\omega$-limit set of $u$.

\item Furthermore, since the space dimension is equal to one, we can use the equivalence between scalar conservation laws and Hamilton-Jacobi equations. More precisely, we can write $u=p + \partial_x U$ for some $U\in L^\infty_{\text{loc}}(W^{1,\infty}(\T))$. Then $U$ is the solution of the Hamilton-Jacobi equation
\[
\p_t U + \frac{1}{2}(p+ \p_x U)^2=\cos^2\left( \frac{\kappa_0 x}{2}\right).
\]
It is well known, see \cite{Ishii, Lions1}, that the Cauchy problem for Hamilton-Jacobi equations with initial data $U_0 \in UC([0, \infty) \times \mathbb{R})$  admits a unique \emph{viscosity solution} $U$.
Investigating the long time behavior of $u$ is therefore equivalent, in some sense, to investigating  the long-time behavior of  $U$.

\item Solutions of Hamilton-Jacobi equations behave, for long times, as ``wave solutions'', which in turn, have a strong connection to solutions of the cell-problem in the homogenization of Hamilton-Jacobi equations.
Hence we will also rely on homogenization results.

\end{itemize}

Let us now give more details on each of the points sketched above.
As a preliminary step, we prove that if $u_0\in L^\infty(\T)$, then $u\in L^\infty(\R_+\times \T)$. First, note that if  $k>0$, then $u^\pm_k(x):=\pm \sqrt{2} (k+ \cos^2(\kappa_0 x/2))^{1/2}$ is a smooth stationary solution of \eqref{Burgers}, and that $\lim_{k\to +\infty} \inf_{x\in \T} u^+_k(x)=+\infty$, $\lim_{k\to + \infty} \sup_{x\in \T} u^-_k(x)=-\infty$. 
Consequently, if $u_0\in L^\infty(\T)$, there exist $k_+, k_->0$ such that $u^-_{k_-}\leq u_0 \leq u^+_{k_+}$. By the maximum principle, this inequality is preserved by the evolution: for all $t\geq 0$, $u^-_{k_-}\leq u(t)\leq u^+_{k_+}$. Thus $u\in L^\infty(\R_+\times \T)$.
Hence it is enough to prove Theorem \ref{thm-Burgers} with $q=1$.

$\bullet$ \textbf{Compactness of the family $(u(t))_{t\in \R}$:}

Let us recall a few results around regularizing effects in scalar conservation laws. In one space dimension and for strictly convex fluxes (which is the case considered here), in the case when $f\equiv 0$, a smoothing effect in $BV$ spaces had been established by Oleinik \cite{Oleinik} and Lax \cite{Lax}.
This property was then generalized to higher space dimensions by Lions, Perthame and Tadmor \cite{LPT}, using the kinetic formulation of the equation and averaging lemmas. 
Bourdarias, Gisclon and Junca \cite{BGJ} and Golse and Perthame \cite{GP} established optimal regularity results in one space dimension, with a regularity index depending on the degeneracy of the flux.
Recently, Gess and Lamy \cite{GL} proved similar results in all dimensions, in the case of a forced scalar conservation law, which is precisely the setting we are considering. 

Since $u\in L^\infty(\R_+\times \T)$, as a  consequence of Theorem 1 of \cite{GL}, we have $u\in W^{s, 3/2}_\text{loc}(\R_+\times \T)$ for all $s<1/3$. Note that in the present case, this result actually follows from the use of the kinetic formulation and from Theorem B in \cite{LPT}. Furthermore,
for any sequence $(t_n)_{n\in \N}$ increasing and converging towards infinity, for any $T>0$, the family $u(t_n+\cdot)$ is uniformly bounded in $W^{s, 3/2}([0,T]\times \T)$. It follows that the sequence $u(t_n + \cdot)$ is compact in $\mathcal C([0,T], L^1(\T))$. Let us denote by $\bar u$ its limit, up to the extraction of a subsequence. Then $\bar u$ is still an entropy solution of \eqref{Burgers}. We will now use the equivalence of \eqref{Burgers} with Hamilton-Jacobi equations to justify that $\bar u$ is in fact a \textit{stationary} solution of \eqref{Burgers}.

$\bullet$ \textbf{Equivalence with Hamilton-Jacobi equations:}

As indicated above, we write $u=p+ \partial_x U$, $\bar u=p+ \partial_x \bar U$. Then $U$ and $\bar U$ are solutions of the Hamilton-Jacobi equation\footnote{Note that $U$ and $\bar U$ are defined up to a function of $t$, and we can always choose this function so that \eqref{HJ} is satisfied.}
\begin{equation}\label{HJ}
\p_t U + H(x, p + \p_x U)=0,
\end{equation}
where the Hamiltonian $H$ is defined by
%
%
%
\begin{equation}\label{def:Hamiltonian}
H(x, q)=\frac{1}{2}q^2-\cos^2 \left(\frac{\kappa_0 x}{2}\right):=T(q)-V(x),\quad x\in \T, \ q\in \R,
\end{equation}
and $T, V$ are respectively the kinetic and the potential energy. 
We already know that up to the extraction of a subsequence, $\p_x U(t_n+\cdot )\to \p_x \bar U$  in $\mathcal C([0,T], L^1(\T))$. On the other hand, the long time behavior of Hamilton-Jacobi equations in the space of Bounded Uniformly Continuous functions $BUC(\T)$ has been investigated by numerous authors, see for instance  \cite{barles2000large,Fathi,namah1999remarks}. 

We will rely in the present note on the following result of Roquejoffre \cite{Roquejoffre}:

\begin{theorem}\label{thm-roquejoffre}
Consider a smooth Hamiltonian $H=H(x, q)$, which is strictly convex and coercive with respect to the variable $p$, i.e. 
\begin{equation}
    \label{coercivity}\lim_{|q| \rightarrow + \infty} \inf_{x\in \mathbb T}H(x, q)=+\infty.
\end{equation}
Let $U_0 \in BUC( \T)$, and let $U(t, x) \in BUC([0, \infty) \times \mathbb{T})$ be the unique viscosity solution to \eqref{HJ} with initial datum $U_0$. Then, there exists a wave solution to \eqref{HJ} of the form
$-\lambda t + \phi (x),$
such that
$$\lim_{t \rightarrow + \infty} \|U(t, \cdot)+\lambda t - \phi\|_{L_x^\infty}=0.$$
\end{theorem}

The above theorem implies that $u=p+ \p_x U $ converges, in the sense of distributions, towards $p+\p_x \phi$, and thus, by uniqueness of the limit in the sense of distributions, $\bar u(t,x)=p+ \p_x \phi(x)$. In particular, $\bar u$ is a stationary entropy solution of the forced Burgers equation \eqref{Burgers}.

Let us now look at the equation satisfied by $\lambda$ and $\phi$. By identification, it is easily proved that $(\lambda, \phi)$ is a solution of
\begin{equation}
\label{eq:cell}
H(x, p + \p_x \phi(x))= \lambda,\quad x\in \T.
\end{equation}
This equation is known as the ``cell-problem'' in homogenization theory. We refer to the seminal paper by Lions, Papanicolaou and Varadhan \cite{Lions}, of which we now recall the main results.

$\bullet$ \textbf{Cell problem and homogenization of Hamilton-Jacobi equations:}

In \cite{Lions}, the authors show that for all $p\in \R$, there exists a unique $\lambda \in \mathbb{R}$ such that there exists $\phi$ viscosity solution to \eqref{eq:cell}. More precisely, they prove the following result.
\begin{theorem}\label{thm-lions1}
Assume that the Hamiltonian $H=H(x,q)$ defined on $\mathbb{T} \times \mathbb{R}$ is periodic in $x$, strictly convex and coercive in $q$.
Then, for each $p \in \mathbb{R}$, there exists a unique $\lambda:=\overline{H}(p) \in \mathbb{R}$, such that there exists $\phi \in \mathcal C(\T)$ a periodic viscosity solution to \eqref{eq:cell}. Moreover, $\overline{H}(p)$ is continuous in $p$.
\end{theorem}

We can explicitly write the computations for the Hamiltonian $H(x, q)$ given by \eqref{def:Hamiltonian}. Following \cite{Lions}, we claim that

\begin{align}
& \label{eq:sol-psmall} \overline{H}(p)=0 \quad \text{if} \quad |p| \le \frac{2\sqrt2}{\pi}=\frac{\sqrt{2}}{2\pi} \int_0^{2\pi} \sqrt{V}(x) \, dx,\\
& \label{eq:sol-pbig}\overline{H}(p)=\lambda \quad \text{where} \quad |p|=\frac{\sqrt{2}}{2\pi} \int_{0}^{2\pi} \sqrt{\cos^2\left(\frac{\kappa_0 x}{2}\right)+\lambda} \, dx,\;  (\lambda \ge 0) \quad \text{if} \; |p| \ge \frac{2\sqrt{2}}{\pi}.
\end{align}

For $|p| \le 2\sqrt{2}/\pi$, for any \emph{minimum} $x_0 \in [0, 2\pi]$ such that $\sqrt{V}(x_0)=0$, following \cite{Lions} we introduce a point $\bar x \in [x_0, x_0 + 2\pi]$ such that

\begin{align*}
\int_{x_0}^{\bar x} \left(\sqrt{2} \left| \cos \left( \frac{\kappa_0 x}{2}\right) \right|-p\right) \: dx= \int^{x_0+2\pi}_{\bar x} \left(\sqrt{2} \left| \cos \left( \frac{\kappa_0 x}{2}\right) \right|+p\right) \: dx,\\
\end{align*}
which gives
\begin{align*}
\int_{x_0}^{\bar x} \left|\cos \left( \dfrac{\kappa_0 x}{2}\right) \right|\; dx = 2 + \frac{p\pi}{\sqrt{2}}.
\end{align*}

One  viscosity solution of \eqref{eq:cell} is then the periodic extension of
\begin{equation}\label{eq-sol-p-small}
\phi(x)=\begin{cases}
\int_{x_0}^{ x} \left(\sqrt{2}\sqrt{V}(y)-p\right) \: dy,\quad x_0\leq x \leq \bar x,\\
\int_{x}^{x_0+2\pi} \left(\sqrt{2} \sqrt{V}(y)+p\right) \, dy,\quad \bar x\leq x \leq x_0+2\pi.
\end{cases}
\end{equation}

Notice that $p + \p_x \phi(x)= \sqrt{2 } \sqrt{V}(x)$ for $x\in (x_0, \bar x)$, and 
 $p + \p_x \phi(x)= - \sqrt{2 }\sqrt{V}(x)$ for $x\in (\bar x, x_0 + 2\pi)$. In particular, $\p_x \phi$ has a jump at $x=\bar x$ except for specific values of $p$ for which $V(\bar x)=0.$ We also recall that the solutions of \eqref{eq:cell} are not unique, as emphasized in Remark \ref{rmk-counting-solutions} below and in \cite{Lions,MRTT}. Indeed, to each minimum $x_0$  of $\sqrt{V}(x)$ corresponds (at least) one viscosity solution. We will comment more on the non-uniqueness of solutions of \eqref{eq:cell} in Remark \ref{rmk-counting-solutions} below.

%
In the complementing case $p \ge \frac{2\sqrt{2}}{\pi}$ (assume $p\ge 0$, the case $p \le 0$ being symmetric), take $\lambda \ge 0$  such that $\displaystyle p=\frac{\sqrt{2}}{2\pi}\int_0^{2\pi} \sqrt{V(x) + \lambda} \, dx$. Then one can find $x_0 \in \left[0, 2\pi\right]$ such that 
\begin{align*}
\sqrt{2\left(V(x_0)+\lambda\right)}=p.
\end{align*}
In this case, a viscosity solution to \eqref{eq:cell} is given by
\begin{equation}\label{eq-visc-sol-p-big}
\phi(x)=\int_{x_0}^x \left(\sqrt{2V(y)+\lambda}-p \right)\; dy
\end{equation}
for $x_0\le x \leq x_0 + 2\pi$.
Again, one considers the periodic extension of $\phi(x)$ defined above. Note that in this case, $\phi \in \mathcal C^1(\T)$. 

\begin{remark}[About the (non)-uniqueness of solutions of \eqref{eq:cell}]\label{rmk-counting-solutions}
Notice that there is a difference between the case $p$ ``small'' ($|p|\le {2\sqrt{2}}/{\pi}$) and $p$ ``big'' ($|p|\ge {2\sqrt{2}}/{\pi}$) in terms of \emph{uniqueness} of the solutions.\\ 
More precisely, for $|p| \le {2\sqrt{2}}/{\pi}$, to \emph{any} $x_0$ such that $\sqrt{V}(x_0)=0$ corresponds a viscosity solution to \eqref{eq:cell}. Therefore, for $p$ ``small'', there are at least as many viscosity solutions to \eqref{eq:cell} as the number of minimum points $x_0 \in \T$ such that $\sqrt{V}(x_0)=0$, up to addition of constants. On the other hand, the viscosity solution $\phi(x)$ to \eqref{eq:cell} is \emph{unique} for $|p|\ge {2\sqrt{2}}/{\pi}$ (again, up to constants).

We can also revisit these results  at the level of the conservation law, following the computations in \cite{MRTT}. Indeed, $\bar u = p + \p_x \phi$ is a stationary entropy solution of the Burgers equation \eqref{Burgers}. As a consequence, $\bar u$ satisfies the following properties:
\begin{itemize}
    \item There exists a constant $\lambda\geq 0$ and a function $\eta\in L^\infty(\T)$ such that $\eta(x)\in \{-1,1\}$ a.e. such that
    \[
    \bar u (x)= \eta (x)\sqrt{2 (V(x) + \lambda)};
    \]
    \item Because of the entropy condition, for all $x \in \T$, $[\bar u]_{|x}:=\bar u(x^+)-\bar u(x^-)\leq 0$;
    
    \item $\langle \bar u\rangle=p. $
    \end{itemize}
Now, since $[\bar u]_{|x}= [\eta]_{|x}\sqrt{2 (V(x) + \lambda)}$, if $[\eta]_{|x_0}>0$ for some $x_0\in \T$, then necessarily $V(x_0) + \lambda=0$. Whence $\lambda=0$ and $x_0$ is a minimum of $V$.

It follows that if $V$ has a unique minimum in $\T$, then the viscosity solution of \eqref{eq:cell} (or equivalently the stationary entropy solution of \eqref{Burgers}) is always unique, up to the addition of a constant. In this case, if $|p|\leq \sqrt{2} \langle \sqrt{V}\rangle$, then $\eta$ has two jumps: one positive jump at the unique point $x_0$ where $V$ reaches its minimum, and one negative jump at the point $\bar x$ introduced above.
However, if $V$ has several distinct minima, then in the case $|p|\leq \langle \sqrt{2V}\rangle$, there are several points where $\eta$ can have a positive jump, and the solutions are no longer unique.
If $|p|\geq \langle \sqrt{2V}\rangle$, then $\eta$ must keep a constant value, and the solution is always unique (and smooth).

In our case, where $\sqrt{V}(x)=\left|\cos({\kappa_0x}/{2})\right|$, there are always at least two minima, and therefore viscosity solutions are not unique when $|p|$ is smaller than the critical value $2\sqrt{2}/\pi$. 
\end{remark}

{\bf $\bullet$ Conclusion:}

Therefore, at this stage we have proved that for all $u_0$, there exists a stationary solution $\bar u $ of the Burgers equation \eqref{Burgers} such that $u(t)\to \bar u$ in $L^1(\T)$. When $p$ is small and $\bar u$ is therefore not uniquely determined,
the numerical simulations in \cite{MRTT} show that the long-time limit solution depends on the initial data.

Let us now look at the Fourier spectrum of $u$ and $\bar u$ in the case when $|p|$ is small.
The Fourier spectrum of each solution $\phi$ of \eqref{eq:cell} decays like $|k|^{-2}$. Indeed, $\phi$ is continuous, but its derivative has a jump at  $x=\bar x$. Therefore the spectrum of $\phi$ is similar to the one of the function $|\cdot|$ on $[-\pi, \pi]$ (extended on $\mathbb R$ by periodicity), which is equivalent to $|k|^{-2}$. Similarly, $\bar u$ is discontinuous, and therefore its Fourier spectrum is equivalent to the one of $\mathrm{sgn}(x)$ and  behaves like $|k|^{-1}$. Let us discuss the time dynamics of the spectrum, starting from an initial data $U_0$ which is smooth, and therefore has a strongly decaying Fourier spectrum. It is well-known that the solution $u$ develops shocks in finite time. Hence there exists $T$ such that for $t<T$, the Fourier spectrum $\hat u(t,k)$ has a strong decay, and for $t\geq T$, the Fourier spectrum decays like $|k|^{-1}$.

\section{Resonance phenomena: the case of unsteady forcing}

\label{sec:resonance}

Let us now look at the Burgers equation in the form \eqref{Burgers-bis}, with a forcing $f(x-\omega t)$, with $f$ a smooth periodic function with zero average. 
Up to now, we only considered the case $\omega=0$. 
But actually, this case corresponds to a resonance: indeed, 
the dispersion relation is linear ($\tau=p k$). 
Since the spatial average of $v$ solution of \eqref{Burgers-bis} is zero, its associated time frequency is also zero, and a resonant forcing is steady in time.
Therefore, following \cite{MRTT} we investigate near-resonances, i.e. we assume that $0<\omega$ and $\omega $ is small. 
We seek for a traveling wave solution $v(t,x)=G(z)$ where $z=x-\omega t$. 
This gives 
\begin{equation}
    \label{eq:cell-bis}
    \dfrac{d}{dz}(\dfrac{1}{2}(G(z)+p-\omega)^2)=f(z).
\end{equation}
    As before (see in particular the computations in Remark \ref{rmk-counting-solutions}), it follows that there exists a function $\eta\in L^\infty(\T)$ with values in $\{-1,1\}$ and a number $\lambda\geq 0$ such that
    \[
    G(z)=\omega-p + \eta(z) \sqrt{2V(z)+ \lambda}
    \]
where $V$ is the unique primitive of $f$ such that $\min_{\T} V=0$. The argument is exactly the same as in the previous section: equation \eqref{eq:cell-bis} always has at least one solution. This solution is unique if $|\omega-p|$ is larger than the critical value $\omega_{cr}=\dfrac{1}{2\pi} \int_0^{2\pi} \sqrt{2 V (z)} \, dz$ or if $V$ has a unique minimum over $\T$. It is not unique if $|\omega-p| <\omega_{cr} $ and $V$ has at least two distinct minima (and is not constant).

Therefore, for $|\omega-p| \geq \omega_{cr} $ \emph{there is a rather abrupt change in behavior which, following \cite{MRTT}, we interpret as the boundary of resonance. That is, a sharp transition from resonant behavior (with the forcing continuously pumping energy into the system, which is dissipated by a shock) to non-resonant behavior (with no work done by the forcing) occurs at $\omega=p \pm \omega_{cr}$}.

Let us reinterpret this in Fourier, assuming  that $f$ is supported by one Fourier mode to fix ideas:
\begin{itemize}
    \item Assume that $\omega$ is close to resonance, i.e. $|\omega-p| < \omega_{cr}$. Then $G$ has infinitely many non-zero modes, and $\hat G(k)$ decays like $|k|^{-1}$. Therefore there is a transfer of energy into high frequencies.
    
    \item Assume now that $\omega$ is non-resonant, i.e. $|\omega-p| > \omega_{cr}$. Then $G$ is smooth and $\hat G(k)$ is rapidly decaying in $k$.
    
\end{itemize}

\section{Links between homogenization, long-time behavior and quasi-resonances}

Starting with the homogeneization problem, we briefly review the following result due to Lions, Papanicolaou and Varadhan (see \cite{Lions}).

\begin{theorem}\label{thm-convergence}
	For any initial datum $U_0 \in BUC(\mathbb{T})$ (Bounded Uniformly Continuous functions), the solution $U^\varepsilon$ to 
	\begin{align}\label{eq:homog-problem}
	\p_t U^\varepsilon + H \left( \frac{x}{\varepsilon}, \p_x U^\varepsilon\right)=0,\quad U^\varepsilon_{|t=0}=U_0
	\end{align}
converges, for $\varepsilon \rightarrow 0$ uniformly in space on $\mathbb{T} \times [0, T]$, for any $T<\infty$, towards the unique viscosity solution in $BUC(\mathbb{T} \times [0, T])$ to the following system:
	\begin{equation}\label{eq_limit-HJ}
	\begin{cases}
	\partial_t \overline{U} + \overline{H}(\partial_x \overline{U})=0, \\
	\overline{U}(0,x)=U_0(x).
	\end{cases}
	\end{equation}
	
\end{theorem}

The proof developed in \cite{Lions} makes rigorous the following reasoning.\\
$\bullet$ Consider the solution $\bar U$ to system \eqref{eq_limit-HJ}, whose existence and uniqueness is due to \cite{Lions1, Ishii}. \\
$\bullet$ For each $\partial_x \overline{U}$ fixed, find $\lambda \in \mathbb{R}$ (provided by Theorem \ref{thm-lions1}), such that there exists a viscosity solution $V(x, y)$ (again, due to Theorem \ref{thm-lions1}) to 
\begin{equation}\label{eq:lambda}
H(y, \partial_x\overline{U} (x)+ \partial_y V(x,y))=\bar H(\partial_x \overline{U}(x)):=\lambda.
\end{equation}
At this point, define
$$U^\varepsilon(t, x)=\overline{U}(t, x)+\varepsilon V \left(x, \frac{x}{\varepsilon}\right)$$
where $\overline{U}, V$ solve \eqref{eq_limit-HJ}-\eqref{eq:lambda} respectively.
Now plug $U^\varepsilon$ in the homogenized system \eqref{eq:homog-problem}. Denoting by $y=x/\varepsilon$ the fast variable, we obtain
\begin{align*}
\partial_t U^\varepsilon+H(y, \partial_x U^\varepsilon)&=\partial_t \overline{U}+H\left(y, \partial_x\overline{U}+\partial_y V\right)+O(\varepsilon)\\
&=- \overline{H}(\partial_x \overline{U})+\overline{H}(\partial_x \overline{U})+O(\varepsilon)=O(\varepsilon).
\end{align*}
Starting from the particular case of affine initial data as
\begin{align*}
&U_0:=\alpha+ p|x|
=\begin{cases}
\alpha+ p x, \quad x  \in [0, \pi], \\
\alpha - p x, \quad x \in [-\pi, 0),
\end{cases}\\
& \text{ with periodic extension     }U_0(x+2\pi)=U_0(x),
\end{align*}
then the solution to \eqref{eq_limit-HJ} is given by
%
%
\begin{equation}\label{eq-sol-affine-data}
\overline{U}(t, x)=\alpha + p |x| - t \overline{H}(p).
\end{equation}
The explicit form of solution \eqref{eq-sol-affine-data} is strictly related to the fact that the initial condition $U_0(x)$ is affine, i.e. $\partial_x U_0$ is piecewise constant, and therefore $\overline{H}(\partial_x U_0)$ does not depend on $x$. However, this is enough to prove the theorem, as it can be extended to the case of general initial data belonging to $BUC(\mathbb{T})$, see \cite{Lions}.\\
This way, Theorem \ref{thm-convergence} tells us that the solution $U^\varepsilon$ to the homogenized equation \eqref{HJ} converges, when $\varepsilon \rightarrow 0$, towards the so-called \emph{wave solution} $\overline{U}(t, x)=\alpha + p |x| - t \overline{H}(p)$ which also appears in the long-time behavior of \eqref{HJ}, see Theorem \ref{thm-roquejoffre}. Hence there are strong links between the homogenization and the long time behavior of Hamilton-Jacobi equations. 

$\bullet$ \textbf{Link between the homogenized problem and the long-time behavior}.\\
This fact can be also seen from the following heuristics.
We write explicitly equation \eqref{eq:homog-problem}, where the Hamiltonian is \eqref{def:Hamiltonian},

\begin{align*}
    \partial_t U^\varepsilon + \dfrac{1}{2} (p+\partial_x U^\varepsilon)^2=V(x/\varepsilon).
\end{align*}

We take the derivative in $x$, which gives

\begin{align*}
    \partial_{tx}U^\varepsilon + (p+\partial_x U^\varepsilon) \partial_{xx} U^\varepsilon=\varepsilon^{-1} f(x/\varepsilon).
\end{align*}

The Burgers variable is then $u^\varepsilon=p+\partial_x U^\varepsilon$. The equation reads

\begin{align*}
    \partial_t u^\varepsilon + u^\varepsilon \partial_x u^\varepsilon=\varepsilon^{-1} f(x/\varepsilon).
\end{align*}

Now set $y:=x/\varepsilon$. This yields

\begin{align*}
    \varepsilon \partial_t u^\varepsilon + u^\varepsilon \partial_y u^\varepsilon=f(y).
\end{align*}

Setting also $\tau:=t/\varepsilon$, we have

\begin{align*}
    \partial_\tau u^\varepsilon + u^\varepsilon \partial_y u^\varepsilon=f(y),
\end{align*}

and notice that $\tau \rightarrow \infty$ as $\varepsilon$ vanishes, then the asymptotics $\varepsilon \rightarrow 0$ of  \eqref{eq:homog-problem} can be viewed as an investigation on long times.\\
$\bullet$ \textbf{Similarity and discrepancy between the near-resonance mechanism and the long-time behavior}\\
As already remarked in Section \ref{sec:resonance}, a steady forcing is resonant for the Burgers equation \eqref{Burgers}.
 Then a near-resonant one evolves slowly in time, like
\begin{align}\label{eq:quasi-steady}
    \partial_t u^\varepsilon + \dfrac{1}{2}\partial_x (u^\varepsilon)^2 = \varepsilon^2 f (x, \varepsilon t).
\end{align}
Now set $\tau:=\varepsilon t$, scale $u^\varepsilon$ as $u^\varepsilon=\varepsilon\tilde{u}$ and observe that the previous equation turns into 
\begin{align}\label{eq:time-dep-forcing}
    \partial_\tau \tilde{u}+\tilde{u}\partial_x \tilde{u}=f(x, \tau),
\end{align}
where, again, $t=\tau/\varepsilon \rightarrow \infty$ as $\varepsilon \rightarrow 0$. Therefore the quasi-steady or near-resonant forcing problem \eqref{eq:quasi-steady} can be also seen as a long-time asymptotics (this is the ``similarity'' part and the motivation to rely on the study of the long-time/homogeneized problem).\\
On the other hand, the fact that the quasi-steady equation \eqref{eq:quasi-steady} is equivalent to \eqref{eq:time-dep-forcing} with $\varepsilon=1$ strongly indicates that the quasi-resonant mechanism described in Section \ref{sec:resonance}, where the cooperation of dissipation due to shocks and nonlinear transfer acts under a certain threshold $\omega_{cr}$ of the time-frequency of the forcing (or, in other words, of the average of the solution), cannot be seen in general as a limit of an equation with a small parameter, because of possible scaling invariances. This is indeed pointed out in \cite{MRTT}, where the authors claim exactly that \emph{near-resonances in the Burgers equation cannot be defined as an asymptotic limit involving a small parameter $\varepsilon$,} in the sense that we cannot shrink the frequencies by means of a small parameter turning the non-resonant to the resonant ones.
This is indeed due to the fact that the distiction between resonant and non-resonant solutions arises from a finite bifurcation in the behavior of the solutions (as we explicitly see in \eqref{eq:sol-pbig}-\eqref{eq-sol-p-small}).

\section*{Acknowledgement}
This project has received funding from the European Research Council (ERC) under the European Union's Horizon 2020 research and innovation program Grant agreement No 637653, project BLOC ``Mathematical Study of Boundary Layers in Oceanic Motion’’. This work was  supported by the SingFlows project, grant ANR-18-CE40-0027 of the French National Research Agency (ANR).

%
%
%

\end{document}